\documentclass[a4paper,10.5pt]{elsarticle}
\usepackage{geometry}
\usepackage{makecell}
\usepackage{extarrows}
\usepackage{subfigure}
\usepackage{mathrsfs}
\usepackage{lineno}
\usepackage{ulem}
\usepackage{graphicx}
\usepackage{subfigure}
\usepackage{makecell}
\usepackage{fancyhdr}
\usepackage{amssymb}
\usepackage{multicol}
\usepackage{amsmath}
\usepackage[figuresright]{rotating}

\usepackage{hyperref}
\usepackage{cleveref}

\newcommand{\I}{\mathcal{I}}

\usepackage{amsthm}
\theoremstyle{definition}

\newtheorem{thm}{Theorem}[section] 
\newtheorem{defn}[thm]{Definition} 
\newtheorem{lem}[thm]{Lemma}
\newtheorem{exam}[thm]{Example}

\begin{document}
\begin{frontmatter}

\title{Interval-valued fuzzy soft $\beta$-covering approximation spaces}

\author{Shizhan Lu$^{1}$}

\address{\small  $^1$School of Management, Jiangsu University, Zhengjiang 212013, China\\




}


\begin{abstract}

The concept of interval-valued fuzzy soft $\beta$-covering approximation spaces (IFS$\beta$CASs) is introduced to combine the theories of soft sets, rough sets and interval-valued fuzzy sets, and some fundamental propositions concerning interval-valued fuzzy soft $\beta$-neighborhoods and soft $\beta$-neighborhoods of IFS$\beta$CASs are explored. And then four kinds of interval-valued fuzzy soft $\beta$-coverings based fuzzy rough sets are researched. Finally, the relationships of four kinds of interval-valued fuzzy soft $\beta$-coverings based fuzzy rough sets are investigated.

\end{abstract}
\end{frontmatter}

\begin{multicols}{2}
\section{Introduction}

The fuzzy set, as proposed by Zadeh \cite{zadeh1965fuzzy}, stands as a renowned instrument for handling uncertainty with diverse applications \cite{ccalik2021novel,2023Intention,2023Evaluation,lin2021risk,liu2021novel,zhang2024event}. 
Subsequently, Gorzalczany \cite{gorzalczany1987method} introduced the notion of interval-valued fuzzy sets, where the membership degree of set elements lies within the interval $[0,1]$. 
Interval-valued fuzzy sets are adept at handling scenarios where precise probabilities of set membership are elusive, offering instead an interval within which such probabilities are constrained \cite{gorzalczany1987method,pkekala2021inclusion}.

Molodtsov \cite{MOLOD} introduced soft sets as a solution to the challenge of uncertainty. 
Unlike traditional mathematical approaches to uncertainty, soft set theory offers a unique advantage in parameter handling. 
Soft set theory enjoys wide-ranging applications across decision-making, rules mining, machine learning, artificial intelligence, image processing, and beyond \cite{zahedi,yang2013,feng2016soft,albaity2023impact,movckovr2021approximations}.

Pawlak \cite{PALA82} initially proposed rough set theory to address the challenges posed by vagueness and granularity in information systems and data analysis. 
Lower and upper approximations are central to rough set theory, offering a mechanism to represent uncertain knowledge based on existing information. 
At its essence, rough set theory constitutes an approximation space defined by a specified universe and an equivalence relation \cite{YUG}. 
Rough set theory has risen to prominence as a rapidly expanding academic discipline spanning domains such as machine learning, deep learning, data mining, decision-making, and beyond \cite{BJAT,KES,MTNA,CHLT,WGLT,guo2022three,JYE,ZKZJ}. 
The amalgamation of soft sets and rough sets frequently inspires the exploration of theories related to soft covering-based rough sets \cite{alcantud2019n,atef2021three,cheng2023error,fan2024overlap,zhang2020fuzzy}, attaining substantial relevance in specific domains. 
However, in fuzzy environments, rough set theory demonstrates inherent limitations, as discussed in \cite{zhangzhan}. 
To overcome these challenges, Zhang and Zhan \cite{zhangzhan} integrated fuzzy sets, soft sets, and rough sets, expanding the notion of soft covering to fuzzy soft $\beta$-covering, drawing insights from the works of Ma \cite{ma2016two} and Yang \cite{yang2016fuzzy,yang2017some}.

To combine the theories of soft sets, rough sets and interval-valued fuzzy sets, the concept of interval-valued fuzzy soft $\beta$-covering approximation spaces (IFS$\beta$CASs) is introduced and the propositions about neighborhoods in IFS$\beta$CASs are researched. And then four kinds of interval-valued fuzzy soft $\beta$-coverings based fuzzy rough sets are researched. Finally, the relationships of four kinds of interval-valued fuzzy soft $\beta$-coverings based fuzzy rough sets are investigated.

The following sections of this manuscript are structured as follows: Section 2 offers the Preliminaries, Section 3 introduces the notion of interval-valued fuzzy soft $\beta$-covering approximation spaces (IFS$\beta$CASs) and investigates the propositions regarding interval-valued fuzzy soft $\beta$-neighborhoods and soft $\beta$-neighborhoods of IFS$\beta$CASs, Section 4 presents four kinds of interval-valued fuzzy soft $\beta$-coverings based fuzzy rough sets  and  studies their relationships.

\section{Preliminaries}

In this section we briefly review some foundations.
If two sets $A$ and $B$ are classical sets (set theory established by Cantor \cite{FF}), $A\sqcap B$ and $A\sqcup B$ represent the intersection and union of $A$ and $B$, respectively. Furthermore, $A\sqsubset B$ represents that $A$ is a subset of the classical set $B$.
Let $U$ be a universal set and $E$ be a set of parameters. Let $IF(U)$ be the set of all interval-valued
fuzzy sets defined over $U$.

\subsection{Fuzzy sets, soft sets, rough sets and approximation spaces}

\begin{defn} \cite{zadeh1965fuzzy} Let $U$ be an initial universe set. A fuzzy set $F$ on $U$ is a mapping $F:U\rightarrow [0,1]$.
\end{defn}

\begin{defn} \cite{MOLOD} Let $U$ be an initial universe set and $E$ be a set of parameters.  $P(U)$ is denoted the power set of $U$. A pair $(F,A)$ is called a soft set over $U$, where $A\subset E$ and $F$  is a mapping given by $F: A\rightarrow P(U)$.
\end{defn}

\begin{defn} \cite{PALA82} Let $U$ be an initial universe set and $R$ be an equivalence relation on $U$. A pair $(U, R)$ is called
a Pawlak approximation space. Two rough approximations $R_-(X)=\{x\in U:[x]_R\sqsubset X\}$ and $R_+(X)=\{x\in U:[x]_R\sqcap X\neq\emptyset\}$ are called the lower and upper approximation of $X$, respectively, where $X\sqsubset U$.
\end{defn}

\begin{defn} \cite{ma2016two}  Let $U$ be an initial universe set and $F(U)$ be the set of all fuzzy subsets of $U$. For $\beta\in(0,1]$, $\mathcal{C}=\{C_1,C_2,\cdots,C_m\}$ is called a fuzzy $\beta$-covering of $U$ if $\bigcup_{i=1}^mC(x)\geqslant\beta$ for all $x\in U$, where $\mathcal{C}\sqsubset F(U)$. And $(U,\mathcal{C})$ is called a fuzzy $\beta$-covering approximation space.

If $\beta=1$, then $(U,\mathcal{C})$ is called a fuzzy covering approximation space.
\end{defn}

\subsection{Interval-valued fuzzy sets and interval-valued fuzzy soft sets}

\begin{defn}\cite{bustince2013generation,perez2015ordering} Assuming that $I_1=[I_1^-,I_1^+]$ and $I_2=[I_2^-,I_2^+]$ are two bounded closed interval on the real number field, then $I_1\leqslant I_2$ if and only if $I_1^+\leqslant I_2^+$ and $I_1^-\leqslant I_2^-$.
\end{defn}

\begin{defn}\cite{dai2023interval,pkekala2021inclusion,zhang2009characterization} Assuming that $I_1=[I_1^-,I_1^+]$ and $I_2=[I_2^-,I_2^+]$ are two bounded closed interval on the real number field,  then we call $I_1$ and $I_2$ as two interval values. For any two interval values $I_1$ and $I_2$, the following operations and product order relation are valid:

(1) $I_1\wedge I_2=[I_1^-\wedge I_2^-, I_1^+\wedge I_2^+]$,

(2) $I_1\vee I_2=[I_1^-\vee I_2^-, I_1^+\vee I_2^+]$,

(3) $I_1^c=[1-I_1^+,1-I_1^-]$,

where $x\wedge y=\min\{x,y\}$ and $x\vee y=\max\{x,y\}$ for all $x,y\in[0,1]$.
\end{defn}

\begin{defn}\cite{dai2023interval,gorzalczany1987method,pkekala2021inclusion} An interval-valued fuzzy set $F$ on an universe $U$ is a mapping such that $F:x\rightarrow I$, where $x\in U$ and $I$ is a closed sub-intervals of $[0,1]$.

For any two interval-valued fuzzy sets $F_1$ and $F_2$ on $U$,  the following basic operations are valid  for all $x\in U$:

(1) $(F_1\cap F_2)(x)=F_1(x)\wedge F_2(x)$,

(2) $(F_1\cup F_2)(x)=F_1(x)\vee F_2(x)$,

(3) $F_1^c(x)=[1-F_1^+(x),1-F_1^-(x)]$,

(4) $F_1\subset F_2\Leftrightarrow F_1(x)\leqslant F_2(x)$.
\end{defn}

\begin{defn} \cite{feng2010application,jiang2013entropy,peng2018algorithms}  Let $U$ be an initial universe set and $E$ be a set of parameters. A pair $(F,E)$ is an interval-valued fuzzy soft sets if $F(e)\in IF(U)$  for every $e\in E$.
\end{defn}

For two interval-valued fuzzy soft sets $(F,A)$ and $(G,B)$ over a common universe $U$, $A\sqsubset B\sqsubset E$. We denote $(F,A)\subset (G,B)$ if $F(e)\subset G(e)$ for every $e\in A$.

\begin{thm}\label{thm8} \cite{bayramov2023interval,mondal1999topology} Let $A,B$ and $C$ be three interval-valued fuzzy soft sets, the following statements hold.

(1) $A\cap B= B\cap A$, $A\cup B= B\cup A$.

(2) $(A\cup B)\cap C=(A\cap C)\cup(B\cap C)$, $(A\cap B)\cup C=(A\cup C)\cap(B\cup C)$.

(3) $(A\cap B)\cap C= A\cap(B\cap C)$, $(A\cup B)\cup C= A\cup(B\cup C)$.

(4) $(A^c)^c= A$.

(5) $(A\cap B)^c= A^c\cup B^c$, $(A\cup B)^c= A^c\cap B^c$.

\end{thm}

\section{Neighborhoods in IFS$\beta$CASs}

\subsection{Interval-valued fuzzy soft $\beta$-neighborhoods in IFS$\beta$CASs}

Let the interval-valued number $\beta$ be an indicator to survey other interval-valued fuzzy sets in the processes of decision-making. $I^U(x)=[1,1]$ and $I^{\emptyset}(x)=[0,0]$ for all $x\in U$.

\begin{defn} Suppose $U$ and $E$ are the sets of objects and parameters, respectively. Let $A\sqsubset E$, $F:E\rightarrow IF(U)$. If $I^U\subset\bigcup\limits_{e\in A}F(e)$, then we called $(F,A)$ an interval-valued fuzzy soft covering over $U$.
\end{defn}

Here, $(U,F,A)$ is called an interval-valued fuzzy soft covering approximation space.

\begin{defn} Let $\beta$ be an interval-valued number, $F:E\rightarrow IF(U)$.
If $\beta\leqslant(\bigcup\limits_{e\in A}F(e))(x)$ for all $x\in U$, then $(F,A)$ is called an interval-valued fuzzy soft $\beta$-covering over $U$, the triple $T=(U,F,A)_{\beta}$ is called an interval-valued fuzzy soft $\beta$-covering approximation space.

\end{defn}

\begin{defn}\label{def5} Let $\beta$ be an interval-valued number, and $U$ and $E$ be the sets of objects and parameters, respectively. Let $e\in A\sqsubset E$, $F:E\rightarrow IF(U)$.
$\widetilde{SN}^{\beta}_x=\bigcap\{F(e)\in F(A):\beta\leqslant F(e)(x)\}$ is called an interval-valued fuzzy soft $\beta$-neighborhood of $x$ in  $(U,F,A)_{\beta}$.

\end{defn}

\begin{lem}\label{lem1} Let $I_*$ be an  interval-valued number, $\I$ and $\I_1$ be two families of interval-valued numbers, $|\I|=\Lambda$, $|\I_1|=\Gamma$, and $\I\sqsubset \I_1$, where $\Lambda$ and $\Gamma$ are two random sets of indicators.

(1) $I_*\leqslant I$ for all $I\in\I$ if and only if $I_*\leqslant\bigwedge\{I: I\in \I\}$.

(2) If $I_*\leqslant I_{\alpha}$ for an $I_{\alpha}\in\I$, then $I_*\leqslant\bigvee\{I: I\in \I\}$.

(3) If $\I\sqsubset\I_1$, then $\bigwedge\{I: I\in \I_1\}\leqslant\bigwedge\{I: I\in \I\}$.

\end{lem}

{\bf\slshape Proof} (1) $I_*\leqslant I$ for all $I\in\I$, i.e., $I_*^-\leqslant I^-$ and $I_*^+\leqslant I^+$ for all $I\in\I$.

$I_*^-\leqslant inf\{I^-:I\in\I\}$ if and only if $I_*^-\leqslant I^-$ for all $I\in\I$. $I_*^+\leqslant inf\{I^+:I\in\I\}$ if and only if $I_*^+\leqslant I^+$ for all $I\in\I$. That is, $I_*\leqslant\bigwedge\{I: I\in \I\}$ if and only if $I_*\leqslant I$ for all $I\in\I$.

(2) $I_*\leqslant I_{\alpha}$ for an $I_{\alpha}\in\I$, i.e., $I_*^-\leqslant I_{\alpha}^-$ and $I_*^+\leqslant I_{\alpha}^+$. Then, $I_*^-\leqslant I_{\alpha}^-\leqslant sup\{I^-:I\in\I\}$ and $I_*^+\leqslant I_{\alpha}^+\leqslant sup\{I^+:I\in\I\}$, it means that  $I_*\leqslant\bigvee\{I: I\in \I\}$.

(3) $\I\sqsubset\I_1$, then $inf\{I^-: I\in \I_1\}\leqslant inf\{I^-: I\in \I\}$ and $inf\{I^+: I\in \I_1\}\leqslant inf\{I^+: I\in \I\}$. It means that $\bigwedge\{I: I\in \I_1\}\leqslant\bigwedge\{I: I\in \I\}$.$\blacksquare$

\begin{thm}\label{pro7} Let $\beta$ be an interval-valued number, $\beta\leqslant\widetilde{SN}^{\beta}_x(x)$ for all $x\in U$
in $(U,F,A)_{\beta}$.

\end{thm}

{\bf\slshape Proof} For each $x\in U$, $\widetilde{SN}^{\beta}_x=\bigcap\{F(e)\in F(A):\beta\leqslant F(e)(x)\}$. By Lemma \ref{lem1} (1), if $\beta\leqslant F(e)(x)$ for all $e\in A^*\sqsubset A$, then $\beta\leqslant (\bigcap\limits_{e\in A^*}F(e))(x)$, i.e., $\beta\leqslant\widetilde{SN}^{\beta}_x(x)$.$\blacksquare$

\begin{thm}\label{prop1} Let $\beta$ be an interval-valued number.
If $\beta\leqslant\widetilde{SN}^{\beta}_x(y)$ and $\beta\leqslant\widetilde{SN}^{\beta}_y(z)$, then $\beta\leqslant\widetilde{SN}^{\beta}_x(z)$ for all $x,y,z\in U$.

\end{thm}

{\bf\slshape Proof}  For each $e\in A\in (U,F,A)_{\beta}$, $\beta\leqslant\widetilde{SN}^{\beta}_x(y)=(\bigcap\limits_{\beta\leqslant F(e)(x)}F(e))(y)$.
By the  sufficiency of Lemma \ref{lem1} (1), we can obtain that $\beta\leqslant F(e)(y)$ for each $F(e)\in\{F(e):\beta\leqslant F(e)(x)\}$. It means that if $\beta\leqslant F(e)(x)$  then $\beta\leqslant F(e)(y)$  for each $e\in A$.

In the similar way above, $\beta\leqslant\widetilde{SN}^{\beta}_y(z)$ implies that if $\beta\leqslant F(e)(y)$  then $\beta\leqslant F(e)(z)$ for each $e\in A$.

To sum up, $\beta\leqslant\widetilde{SN}^{\beta}_x(y)$ and $\beta\leqslant\widetilde{SN}^{\beta}_y(z)$ imply that if $\beta\leqslant F(e)(x)$  then $\beta\leqslant F(e)(z)$  for each $e\in A$.
By the necessity of Lemma \ref{lem1} (1), $\beta\leqslant(\bigcap\{F(e):\beta\leqslant F(e)(x)\})(z)=\widetilde{SN}^{\beta}_x(z)$.$\blacksquare$

\begin{thm}\label{prop2} For two interval-valued numbers $\beta_1$ and $\beta_2$,
if $\beta_1\leqslant\beta_2$, then $\widetilde{SN}^{\beta_1}_x\subset\widetilde{SN}^{\beta_2}_x$ for all $x\in U$.

\end{thm}

{\bf\slshape Proof}   Since $\beta_1\leqslant \beta_2$, if $\beta_2\leqslant F(e)(x)$  then  $\beta_1(x)\leqslant F(e)(x)$ for  $e\in A\in (U,F,A)_{\beta}$, i.e., $\{F(e):\beta_2\leqslant F(e)(x),e\in A\}\sqsubset\{F(e):\beta_1\leqslant F(e)(x),e\in A\}$.

By Lemma \ref{lem1} (3), $\{F(e):\beta_2\leqslant F(e)(x),e\in A\}\sqsubset\{F(e):\beta_1\leqslant F(e)(x),e\in A\}$ implies $\widetilde{SN}^{\beta_1}_x=\bigcap\{F(e):\beta_1\leqslant F(e)(x),e\in A\}\subset\bigcap\{F(e):\beta_2\leqslant F(e)(x),e\in A\}=\widetilde{SN}^{\beta_2}_x$.$\blacksquare$

\begin{thm}\label{pro8} Let $\beta$ be an interval-valued number,  $x,y\in U$.

(1)  $\beta\leqslant \widetilde{SN}^{\beta}_x(y)$ if and only if $\widetilde{SN}^{\beta}_y\subset \widetilde{SN}^{\beta}_x$.

(2) $\beta\leqslant \widetilde{SN}^{\beta}_x(y)$ and $\beta\leqslant \widetilde{SN}^{\beta}_y(x)$ if and only if $\widetilde{SN}^{\beta}_x=\widetilde{SN}^{\beta}_y$.

\end{thm}

{\bf\slshape Proof} (1) For $e\in A\in (U,F,A)_{\beta}$, if $\beta\leqslant \widetilde{SN}^{\beta}_x(y)=(\bigcap\limits_{\beta\leqslant F(e)(x)}F(e))(y)$, by the sufficiency of Lemma \ref{lem1} (1), $\beta\leqslant F(e)(y)$ holds for each $F(e)\in\{F(e):\beta\leqslant F(e)(x)\}$, then   $\{F(e):\beta\leqslant F(e)(x)\}\sqsubset\{F(e):\beta\leqslant F(e)(y)\}$ is obtained.

By Lemma \ref{lem1} (1), $\{F(e):\beta\leqslant F(e)(x)\}\sqsubset\{F(e):\beta\leqslant F(e)(y)\}$ implies $\bigcap\{F(e):\beta\leqslant F(e)(y)\}\subset \bigcap\{F(e):\beta\leqslant F(e)(x)\}$, i.e., $\widetilde{SN}^{\beta}_y\subset \widetilde{SN}^{\beta}_x$.

On the other hand, by Theorem \ref{pro7} (1), $\beta\leqslant\widetilde{SN}^{\beta}_y(y)$ holds. If $\widetilde{SN}^{\beta}_y\subset \widetilde{SN}^{\beta}_x$, then $\widetilde{SN}^{\beta}_y(y)\leqslant \widetilde{SN}^{\beta}_x(y)$, then  $\beta\leqslant\widetilde{SN}^{\beta}_y(y)\leqslant\widetilde{SN}^{\beta}_x(y)$.

(2) By the results of (1), $\beta\leqslant \widetilde{SN}^{\beta}_x(y)$ and $\beta\leqslant \widetilde{SN}^{\beta}_y(x)$ if and only if $\widetilde{SN}^{\beta}_y\subset \widetilde{SN}^{\beta}_x$ and $\widetilde{SN}^{\beta}_x\subset \widetilde{SN}^{\beta}_y$, i.e., $\widetilde{SN}^{\beta}_y=\widetilde{SN}^{\beta}_x$.$\blacksquare$

\subsection{Soft $\beta$-neighborhoods in IFS$\beta$CASs}

\begin{defn}\label{def6} Let $\beta$ be an interval-valued number,  $F:E\rightarrow IF(U)$, $A\sqsubset E$.
$\overline{SN}^{\beta}_x=\{y\in U: \beta\leqslant \widetilde{SN}^{\beta}_x(y)\}$ is called the  soft $\beta$-neighborhood of $x$ in $(U,F,A)_{\beta}$.
\end{defn}

\begin{thm}\label{pro320} Let $\beta$ be an interval-valued number,
$x\in\overline{SN}^{\beta}_x$ for all $x\in U$.

{\bf\slshape Proof}  By Theorem \ref{pro7}, $\beta\leqslant\widetilde{SN}^{\beta}_x(x)$, then $x\in\overline{SN}^{\beta}_x$.$\blacksquare$

\end{thm}

\begin{thm}\label{pro321} Let $\beta$ be an interval-valued number,
 $y\in\overline{SN}^{\beta}_x$ if and only if $\overline{SN}^{\beta}_y\sqsubset\overline{SN}^{\beta}_x$ for all  $x,y\in U$.

{\bf\slshape Proof}
If $y\in\overline{SN}^{\beta}_x$, i.e., $\beta\leqslant \widetilde{SN}^{\beta}_x(y)$,  by Theorem \ref{pro8} (1) we have  $\widetilde{SN}^{\beta}_y\subset \widetilde{SN}^{\beta}_x$, then $\widetilde{SN}^{\beta}_y(z)\leqslant \widetilde{SN}^{\beta}_x(z)$ for each $z\in U$.

If $\beta\leqslant\widetilde{SN}^{\beta}_y(z)$, then $\beta\leqslant\widetilde{SN}^{\beta}_y(z)\leqslant\widetilde{SN}^{\beta}_x(z)$. It means that if $z\in\overline{SN}^{\beta}_y$ then $z\in\overline{SN}^{\beta}_x$, i.e., $\overline{SN}^{\beta}_y\sqsubset\overline{SN}^{\beta}_x$.

On the other hand, by the results of Theorem \ref{pro320}, $y\in\overline{SN}^{\beta}_y$. If $\overline{SN}^{\beta}_y\sqsubset\overline{SN}^{\beta}_x$, then $y\in\overline{SN}^{\beta}_x$.$\blacksquare$

\end{thm}

\begin{thm}\label{pro322} Let $\beta$ be an interval-valued number,
if $y\in\overline{SN}^{\beta}_x$ and $z\in\overline{SN}^{\beta}_y$, then $z\in\overline{SN}^{\beta}_x$ for all $x,y,z\in U$.

{\bf\slshape Proof}  $y\in\overline{SN}^{\beta}_x$ and $z\in\overline{SN}^{\beta}_y$, i.e., $\beta\leqslant \widetilde{SN}^{\beta}_x(y)$ and $\beta\leqslant \widetilde{SN}^{\beta}_y(z)$. By Theorem \ref{prop1}, if $\beta\leqslant\widetilde{SN}^{\beta}_x(y)$ and $\beta\leqslant\widetilde{SN}^{\beta}_y(z)$, then $\beta\leqslant\widetilde{SN}^{\beta}_x(z)$, i.e., $z\in\overline{SN}^{\beta}_x$.$\blacksquare$

\end{thm}

\begin{thm}\label{pro322n1} Let $\beta$ be an interval-valued number,  the following statements hold for all $x,y\in U$.

(1)  $\widetilde{SN}^{\beta}_y\subset \widetilde{SN}^{\beta}_x$ if and only if $\overline{SN}^{\beta}_y\sqsubset\overline{SN}^{\beta}_x$.

(2)  $\overline{SN}^{\beta}_x=\overline{SN}^{\beta}_y$ if and only if $\widetilde{SN}^{\beta}_x= \widetilde{SN}^{\beta}_y$.

{\bf\slshape Proof} (1) By Theorem \ref{pro8} (1), $\widetilde{SN}^{\beta}_y\subset \widetilde{SN}^{\beta}_x\Leftrightarrow\beta\leqslant \widetilde{SN}^{\beta}_x(y)$. By Definition \ref{def6}, $\beta\leqslant \widetilde{SN}^{\beta}_x(y)\Leftrightarrow y\in\overline{SN}^{\beta}_x$. By Theorem \ref{pro321}, $y\in\overline{SN}^{\beta}_x\Leftrightarrow\overline{SN}^{\beta}_y\sqsubset\overline{SN}^{\beta}_x$.
Then,  $\widetilde{SN}^{\beta}_y\subset \widetilde{SN}^{\beta}_x\Leftrightarrow\overline{SN}^{\beta}_y\sqsubset\overline{SN}^{\beta}_x$.

(2) By the results of (1), $\overline{SN}^{\beta}_y\sqsubset\overline{SN}^{\beta}_x$ and $\overline{SN}^{\beta}_x\sqsubset\overline{SN}^{\beta}_y$ if and only if $\widetilde{SN}^{\beta}_y\subset \widetilde{SN}^{\beta}_x$ and $\widetilde{SN}^{\beta}_x\subset \widetilde{SN}^{\beta}_y$. It means that $\overline{SN}^{\beta}_y=\overline{SN}^{\beta}_x$ if and only if $\widetilde{SN}^{\beta}_x=\widetilde{SN}^{\beta}_y$.$\blacksquare$

\end{thm}

\begin{thm}\label{pro322n} Let $\beta$ be an interval-valued number, $\Lambda$ is a random set of indicators,   the following statements hold for all $x,y\in U$.

(1)  $\overline{SN}^{\beta}_x\sqcup\overline{SN}^{\beta}_y\sqsubset\overline{\widetilde{SN}^{\beta}_x\cup\widetilde{SN}^{\beta}_y}^{\beta}$.

(2) $\overline{SN}^{\beta}_x\sqcap\overline{SN}^{\beta}_y=\overline{\widetilde{SN}^{\beta}_x\cap\widetilde{SN}^{\beta}_y}^{\beta}$.

(3) $\sqcup_{v\in\Lambda}\overline{SN}_{x_v}^{a\beta}\sqsubset\overline{\bigcup_{v\in\Lambda}\widetilde{SN}_{x_v}^{a\beta}}^{a\beta}$.

(4)  $\sqcap_{v\in\Lambda}\overline{SN}_{x_v}^{a\beta}=\overline{\bigcap_{v\in\Lambda}\widetilde{SN}_{x_v}^{a\beta}}^{a\beta}$.

{\bf\slshape Proof} (1) For each $z\in U$, if $z\in\overline{SN}^{\beta}_x\sqcup\overline{SN}^{\beta}_y$, i.e., $\beta\leqslant\widetilde{SN}^{\beta}_x(z)$ or $\beta\leqslant\widetilde{SN}^{\beta}_y(z)$. By Lemma \ref{lem1} (2), $\beta\leqslant(\widetilde{SN}^{\beta}_x\cup\widetilde{SN}^{\beta}_y)(z)$, i.e., $z\in\overline{\widetilde{SN}^{\beta}_x\cup\widetilde{SN}^{\beta}_y}^{\beta}$. Then, $\overline{SN}^{\beta}_x\sqcup\overline{SN}^{\beta}_y\sqsubset\overline{\widetilde{SN}^{\beta}_x\cup\widetilde{SN}^{\beta}_y}^{\beta}$.

(2) By Lemma \ref{lem1} (1), $\beta\leqslant\widetilde{SN}^{\beta}_x(z)$ and $\beta\leqslant\widetilde{SN}^{\beta}_y(z)$ if and only if  $\beta\leqslant(\widetilde{SN}^{\beta}_x\cap\widetilde{SN}^{\beta}_y)(z)$. Then $\overline{\widetilde{SN}^{\beta}_x\cap\widetilde{SN}^{\beta}_y}^{\beta}=\overline{SN}^{\beta}_x\sqcap\overline{SN}^{\beta}_y$.

(3)  For each $z\in U$, if $z\in\sqcup_{v\in\Lambda}\overline{SN}_{x_v}^{\beta}$, then there exists at least one $v'\in\Lambda$ such that $z\in\overline{SN}_{x_{v'}}^{\beta}$, i.e., $\beta\leqslant(\widetilde{SN}_{x_{v'}}^{\beta})(z)$. By Lemma \ref{lem1} (2),
$\beta\leqslant(\bigcup_{v\in\Lambda}\widetilde{SN}_{x_v}^{\beta})(z)$, i.e., $z\in\overline{\bigcup_{v\in\Lambda}\widetilde{SN}_{x_v}^{\beta}}^{\beta}$. Then,  $\sqcup_{v\in\Lambda}\overline{SN}_{x_v}^{\beta}\sqsubset\overline{\bigcup_{v\in\Lambda}\widetilde{SN}_{x_v}^{\beta}}^{\beta}$.

(4) For each $z\in U$,  $z\in\sqcap_{v\in\Lambda}\overline{SN}_{x_v}^{\beta}$ if and only if $z\in\overline{SN}_{x_v}^{\beta}$ for all $v\in\Lambda$, i.e., $\beta\leqslant(\widetilde{SN}_{x_{v}}^{\beta})(z)$ holds for all $v\in\Lambda$.

By Lemma \ref{lem1} (1),  $\beta\leqslant(\bigcap_{v\in\Lambda}\widetilde{SN}_{x_v}^{\beta})(z)$ if and only if $\beta\leqslant(\widetilde{SN}_{x_{v}}^{\beta})(z)$ for all $v\in\Lambda$, i.e., $\sqcap_{v\in\Lambda}\overline{SN}_{x_v}^{\beta}=\overline{\bigcap_{v\in\Lambda}\widetilde{SN}_{x_v}^{\beta}}^{\beta}$.$\blacksquare$

\end{thm}

\begin{exam}  Let $\beta=[0.5,0.6]$. If $\widetilde{SN}_{x}^{\beta}(z)=[0.3,0.6]$ and $\widetilde{SN}_{y}^{\beta}(z)=[0.5,0.55]$, then $(\widetilde{SN}_{x}^{\beta}\cup\widetilde{SN}_{y}^{\beta})(z)=[0.3,0.6]\vee[0.5,0.55]=[0.5,0.6]$. $z\not\in\overline{SN}^{\beta}_x$ and $z\not\in\overline{SN}^{\beta}_y$, however, $z\in\overline{\widetilde{SN}^{\beta}_x\cup\widetilde{SN}^{\beta}_y}^{\beta}$. Hence, the (1) and (3) of Theorem \ref{pro322n}  cannot take the equal sign.

\end{exam}

\section{Four kinds of interval-valued fuzzy soft $\beta$-coverings based fuzzy rough sets and their relationships}

\subsection{The 1st kind of interval-valued fuzzy soft $\beta$-coverings based fuzzy rough sets}

\begin{defn}  Let $\beta$ be an interval-valued number, $X\in IF(U)$, $x,y\in U$. Let
\begin{equation*}
\left\{
\begin{array}{l}
\widetilde{SA}_{-1}^{\beta}(X)(x)=\bigwedge\limits_{y\in U}\{(\widetilde{SN}_x^{\beta})^c(y)\vee X(y)\}, \\
\widetilde{SA}_{+1}^{\beta}(X)(x)=\bigvee\limits_{y\in U}\{\widetilde{SN}_x^{\beta}(y)\wedge X(y)\}.
\end{array}
\right.
\end{equation*}
$\widetilde{SA}_{-1}^{\beta}(X)$ and $\widetilde{SA}_{+1}^{\beta}(X)$ are the lower approximation and upper approximation of the interval-valued fuzzy set $X$ in $(U,F,A)_{\beta}$, respectively.

If $\widetilde{SA}_{-1}^{\beta}(X)\neq\widetilde{SA}_{+1}^{\beta}(X)$, then $X$ is called the 1st type of interval-valued fuzzy
soft $\beta$-covering based interval-valued fuzzy rough sets; otherwise $X$ is called
interval-valued fuzzy definable.

\end{defn}

\begin{thm}\label{proada} Let $\beta$ be an interval-valued number, $X,Y\in IF(U)$, $x,y\in U$,  the following statements hold.

(1)  $\widetilde{SA}_{-1}^{\beta}(I^U)=I^U$, $\widetilde{SA}_{+1}^{\beta}(I^{\emptyset})= I^{\emptyset}$.

(2) $\widetilde{SA}_{-1}^{\beta}(X^c)=(\widetilde{SA}_{+1}^{\beta}(X))^c$, $\widetilde{SA}_{+1}^{\beta}(X^c)=(\widetilde{SA}_{-1}^{\beta}(X))^c$.

(3)  $\widetilde{SA}_{-1}^{\beta}(X\cap Y)=\widetilde{SA}_{-1}^{\beta}(X)\cap \widetilde{SA}_{-1}^{\beta}(Y)$, $\widetilde{SA}_{+1}^{\beta}(X\cup Y)=\widetilde{SA}_{+1}^{\beta}(X)\cup \widetilde{SA}_{+1}^{\beta}(Y)$.

(4) If $X\subset Y$, then $\widetilde{SA}_{-1}^{\beta}(X)\subset\widetilde{SA}_{-1}^{\beta}(Y)$, $\widetilde{SA}_{+1}^{\beta}(X)\subset\widetilde{SA}_{+1}^{\beta}(Y)$.

(5) $\widetilde{SA}_{-1}^{\beta}(X)\cup\widetilde{SA}_{-1}^{\beta}(Y)\subset\widetilde{SA}_{-1}^{\beta}(X\cup Y)$, $\widetilde{SA}_{+1}^{\beta}(X\cap Y)\subset\widetilde{SA}_{+1}^{\beta}(X)\cap\widetilde{SA}_{+1}^{\beta}(Y)$.

(6) If $(\widetilde{SN}_x^{\beta})^c(x)\leqslant X(x)\leqslant\widetilde{SN}_x^{\beta}(x)$ for all $x\in U$, then $\widetilde{SA}_{-1}^{\beta}(X)\subset X\subset\widetilde{SA}_{+1}^{\beta}(X)$.

(7) If $(\widetilde{SN}_x^{\beta})^c(x)\leqslant X(x)\leqslant\widetilde{SN}_x^{\beta}(x)$ for all $x\in U$, then $\widetilde{SA}_{-1}^{\beta}(\widetilde{SA}_{-1}^{\beta}(X))\subset\widetilde{SA}_{-1}^{\beta}(X)\subset X\subset\widetilde{SA}_{+1}^{\beta}(X)\subset\widetilde{SA}_{+1}^{\beta}(\widetilde{SA}_{+1}^{\beta}(X))$.

(8) If $X\subset Y$, $\widetilde{SA}_{-1}^{\beta}(X)\cup\widetilde{SA}_{-1}^{\beta}(Y)=\widetilde{SA}_{-1}^{\beta}(X\cup Y)$, $\widetilde{SA}_{+1}^{\beta}(X\cap Y)=\widetilde{SA}_{+1}^{\beta}(X)\cap\widetilde{SA}_{+1}^{\beta}(Y)$.
\end{thm}

{\bf\slshape Proof}  (1) For all $x,y\in U$, since $I^U(y)=[1,1]$, then $\widetilde{SA}_{-1}^{\beta}(I^U)(x)=\bigwedge\limits_{y\in U}\{(\widetilde{SN}_x^{\beta})^c(y)\vee I^U(y)\}=\bigwedge\limits_{y\in U}\{[1,1]\}=[1,1]=I^U(x)$, i.e., $\widetilde{SA}_{-1}^{\beta}(I^U)=I^U$.

For all $x,y\in U$, since $I^{\emptyset}(y)=[0,0]$, then $\widetilde{SA}_{+1}^{\beta}(I^{\emptyset})(x)=\bigvee\limits_{y\in U}\{\widetilde{SN}_x^{\beta}(y)\wedge I^{\emptyset}(y)\}=\bigvee\limits_{y\in U}\{[0,0]\}=[0,0]=I^{\emptyset}(x)$, i.e., $\widetilde{SA}_{+1}^{\beta}(I^{\emptyset})=I^{\emptyset}$.

(2) We prove the case of $|U|=2$ at first. Suppose $U=\{y_1,y_2\}$.

By the results of Theorem \ref{thm8} (4) and (5), $\widetilde{SA}_{-1}^{\beta}(X^c)(y_1)=\bigwedge\limits_{y_i\in U}\{(\widetilde{SN}_{y_1}^{\beta})^c(y_i)\vee X^c(y_i)\}=[(\widetilde{SN}_{y_1}^{\beta})^c(y_1)\vee X^c(y_1)]\wedge[(\widetilde{SN}_{y_1}^{\beta})^c(y_2)\vee X^c(y_2)]=[\widetilde{SN}_{y_1}^{\beta}(y_1)\wedge X(y_1)]^c\wedge[\widetilde{SN}_{y_1}^{\beta}(y_2)\wedge X(y_2)]^c=\{[\widetilde{SN}_{y_1}^{\beta}(y_1)\wedge X(y_1)]\vee[\widetilde{SN}_{y_1}^{\beta}(y_2)\wedge X(y_2)]\}^c=(\widetilde{SA}_{+1}^{\beta}(X))^c(y_1)$.

$\widetilde{SA}_{-1}^{\beta}(X^c)(y_2)=(\widetilde{SA}_{+1}^{\beta}(X))^c(y_2)$ can be proved as above.

To sum up,  $\widetilde{SA}_{-1}^{\beta}(X^c)=(\widetilde{SA}_{+1}^{\beta}(X))^c$ hold for the case of $|U|=2$. Other cases of $|U|<\infty$ can be proved by using the mathematical induction.

$\widetilde{SA}_{+1}^{\beta}(X^c)(y_1)=\bigvee\limits_{y_i\in U}\{\widetilde{SN}_{y_1}^{\beta}(y_i)\wedge X^c(y_i)\}=[\widetilde{SN}_{y_1}^{\beta}(y_1)\wedge X^c(y_1)]\vee[\widetilde{SN}_{y_1}^{\beta}(y_2)\wedge X^c(y_2)]=\{[(\widetilde{SN}_{y_1}^{\beta})^c(y_1)\vee X(y_1)]\wedge[(\widetilde{SN}_{y_1}^{\beta})^c(y_2)\vee X(y_2)]\}^c=(\widetilde{SA}_{-1}^{\beta}(X))^c(y_1)$.

$\widetilde{SA}_{+1}^{\beta}(X^c)(y_2)=(\widetilde{SA}_{-1}^{\beta}(X))^c(y_2)$ can be proved as above.

To sum up,  $\widetilde{SA}_{+1}^{\beta}(X^c)=(\widetilde{SA}_{-1}^{\beta}(X))^c$ hold for the case of $|U|=2$.
Other cases of $|U|<\infty$ can be proved by using the mathematical induction.

(3) For each $x\in U$, by Theorem \ref{thm8} (2),  $\widetilde{SA}_{-1}^{\beta}(X\cap Y)(x)=\bigwedge\limits_{y\in U}\{(\widetilde{SN}_x^{\beta})^c(y)\vee (X\cap Y)(y)\}=\bigwedge\limits_{y\in U}\{[(\widetilde{SN}_x^{\beta})^c(y)\vee X(y)]\wedge[(\widetilde{SN}_x^{\beta})^c(y)\vee Y(y)]\}=\{\bigwedge\limits_{y\in U}\{(\widetilde{SN}_x^{\beta})^c(y)\vee X(y)\}\}\wedge\{\bigwedge\limits_{y\in U}\{(\widetilde{SN}_x^{\beta})^c(y)\vee Y(y)\}\}=\widetilde{SA}_{-1}^{\beta}(X)(x)\wedge\widetilde{SA}_{-1}^{\beta}(Y)(x)$.  Then,  $\widetilde{SA}_{-1}^{\beta}(X\cap Y)=\widetilde{SA}_{-1}^{\beta}(X)\cap \widetilde{SA}_{-1}^{\beta}(Y)$.

For each $x\in U$,  by Theorem \ref{thm8} (2), $\widetilde{SA}_{+1}^{\beta}(X\cup Y)(x)=\bigvee\limits_{y\in U}\{\widetilde{SN}_x^{\beta}(y)\wedge (X\cup Y)(y)\}=\bigvee\limits_{y\in U}\{[\widetilde{SN}_x^{\beta}(y)\wedge X(y)]\vee[\widetilde{SN}_x^{\beta}(y)\wedge Y(y)]\}=\{\bigvee\limits_{y\in U}\{\widetilde{SN}_x^{\beta}(y)\wedge X(y)\}\}\vee\{\bigcup\limits_{y\in U}\{\widetilde{SN}_x^{\beta}(y)\wedge Y(y)\}\}=\widetilde{SA}_{+1}^{\beta}(X)(x)\vee \widetilde{SA}_{+1}^{\beta}(Y)(x)$. Then, $\widetilde{SA}_{+1}^{\beta}(X\cup Y)=\widetilde{SA}_{+1}^{\beta}(X)\cup \widetilde{SA}_{+1}^{\beta}(Y)$.

(4)  Since $X\subset Y$, thus $(\widetilde{SN}_x^{\beta})^c\cup X\subset (\widetilde{SN}_x^{\beta})^c\cup Y$, then $(\widetilde{SN}_x^{\beta})^c(y)\vee X(y)\leqslant (\widetilde{SN}_x^{\beta})^c(y)\vee Y(y)$ for all $x,y\in U$, it implies that  $\widetilde{SA}_{-1}^{\beta}(X)(x)\leqslant\widetilde{SA}_{-1}^{\beta}(Y)(x)$ for all $x\in U$, i.e., $\widetilde{SA}_{-1}^{\beta}(X)\subset\widetilde{SA}_{-1}^{\beta}(Y)$.

Since  $X\subset Y$, thus $\widetilde{SN}_x^{\beta}\cap X\subset\widetilde{SN}_x^{\beta}\cap Y$, then $\widetilde{SN}_x^{\beta}(y)\wedge X(y)\leqslant\widetilde{SN}_x^{\beta}(y)\wedge Y(y)$ for all $x,y\in U$, it implies that $\widetilde{SA}_{+1}^{\beta}(X)(x)\leqslant \widetilde{SA}_{+1}^{\beta}(Y)(x)$ for all $x\in U$, i.e., $\widetilde{SA}_{+1}^{\beta}(X)\subset\widetilde{SA}_{+1}^{\beta}(Y)$.

(5) $X\subset X\cup Y$ and $Y\subset X\cup Y$. By the results of (4), $\widetilde{SA}_{-1}^{\beta}(X)\subset\widetilde{SA}_{-1}^{\beta}(X\cup Y)$ and $\widetilde{SA}_{-1}^{\beta}(Y)\subset\widetilde{SA}_{-1}^{\beta}(X\cup Y)$. Then  $\widetilde{SA}_{-1}^{\beta}(X)\cup\widetilde{SA}_{-1}^{\beta}(Y)\subset\widetilde{SA}_{-1}^{\beta}(X\cup Y)$.

$X\cap Y\subset X$ and $X\cap Y\subset Y$. By the results of (4), $\widetilde{SA}_{+1}^{\beta}(X\cap Y)\subset\widetilde{SA}_{+1}^{\beta}(X)$ and $\widetilde{SA}_{+1}^{\beta}(X\cap Y)\subset\widetilde{SA}_{+1}^{\beta}(Y)$. Then $\widetilde{SA}_{+1}^{\beta}(X\cap Y)\subset\widetilde{SA}_{+1}^{\beta}(X)\cap\widetilde{SA}_{+1}^{\beta}(Y)$.

(6) Since $(\widetilde{SN}_x^{\beta})^c(x)\leqslant X(x)$ for all $x\in U$, then $(\widetilde{SN}_x^{\beta})^c(x)\vee X(x)\leqslant X(x)\vee X(x)=X(x)$, then  $\widetilde{SA}_{-1}^{\beta}(X)(x)=\bigwedge\limits_{y\in U}\{(\widetilde{SN}_x^{\beta})^c(y)\vee X(y)\}\leqslant(\widetilde{SN}_x^{\beta})^c(x)\vee X(x)\leqslant X(x)$, i.e., $\widetilde{SA}_{-1}^{\beta}(X)\subset X$.

Since $X(x)\leqslant\widetilde{SN}_x^{\beta}(x)$ for all $x\in U$, then $X(x)=X(x)\wedge X(x)\leqslant X(x)\wedge(\widetilde{SN}_x^{\beta})(x)$, then $X(x)\leqslant X(x)\wedge\widetilde{SN}_x^{\beta}(x)\leqslant\bigvee\limits_{y\in U}\{\widetilde{SN}_x^{\beta}(y)\wedge X(y)\}=\widetilde{SA}_{+1}^{\beta}(X)(x)$, i.e., $X\subset\widetilde{SA}_{+1}^{\beta}(X)$.

Hence, $\widetilde{SA}_{-1}^{\beta}(X)\subset X\subset\widetilde{SA}_{+1}^{\beta}(X)$.

(7) $(\widetilde{SN}_x^{\beta})^c(x)\leqslant X(x)\leqslant\widetilde{SN}_x^{\beta}(x)$ for all $x\in U$, by the results of (6), we have $\widetilde{SA}_{-1}^{\beta}(X)\subset X\subset\widetilde{SA}_{+1}^{\beta}(X)$. By the results of (4), we have $\widetilde{SA}_{-1}^{\beta}(\widetilde{SA}_{-1}^{\beta}(X))\subset\widetilde{SA}_{-1}^{\beta}(X)\subset X\subset\widetilde{SA}_{+1}^{\beta}(X)\subset\widetilde{SA}_{+1}^{\beta}(\widetilde{SA}_{+1}^{\beta}(X))$.

(8) If $X\subset Y$, then $X\cup Y\subset Y\cup Y= Y$. By the results of (4), $\widetilde{SA}_{-1}^{\beta}(X\cup Y)\subset\widetilde{SA}_{-1}^{\beta}(Y)\subset\widetilde{SA}_{-1}^{\beta}(X)\cup\widetilde{SA}_{-1}^{\beta}(Y)$.

By the results of (5), $\widetilde{SA}_{-1}^{\beta}(X)\cup\widetilde{SA}_{-1}^{\beta}(Y)\subset\widetilde{SA}_{-1}^{\beta}(X\cup Y)$.

Then, $\widetilde{SA}_{-1}^{\beta}(X)\cup\widetilde{SA}_{-1}^{\beta}(Y)=\widetilde{SA}_{-1}^{\beta}(X\cup Y)$.

If $X\subset Y$, then $X=X\cap X\subset X\cap Y$. By the results of (4), $\widetilde{SA}_{+1}^{\beta}(X)\subset\widetilde{SA}_{+1}^{\beta}(X\cap Y)$. Then, we have $\widetilde{SA}_{+1}^{\beta}(X)\cap\widetilde{SA}_{+1}^{\beta}(Y)\subset\widetilde{SA}_{+1}^{\beta}(X)\subset\widetilde{SA}_{+1}^{\beta}(X\cap Y)$.

By the results of (5), $\widetilde{SA}_{+1}^{\beta}(X\cap Y)\subset\widetilde{SA}_{+1}^{\beta}(X)\cap\widetilde{SA}_{+1}^{\beta}(Y)$.

Hence, $\widetilde{SA}_{+1}^{\beta}(X\cap Y)=\widetilde{SA}_{+1}^{\beta}(X)\cap\widetilde{SA}_{+1}^{\beta}(Y)$.$\blacksquare$

\begin{defn}\label{defcla1}  Let $\beta$ be an interval-valued number and $X\sqsubset U$ be an object set, $x,y\in U$.
$\overline{SA}_{-1}^{\beta}(X)=\{x\in U:\overline{SN}^{\beta}_x\sqsubset X\}$ and $\overline{SA}_{+1}^{\beta}(X)=\{x\in U:\overline{SN}^{\beta}_x\sqcap X\neq\emptyset\}$ are the lower approximation and upper approximation of the object set $X$ in $(U,F,A)_{\beta}$, respectively.

If $\overline{SA}_{-1}^{\beta}(X)\neq \overline{SA}_{+1}^{\beta}(X)$, then $X$ is called the 1st type of soft
$\beta$-covering based rough sets. Otherwise, $X$ is called definable.
\end{defn}

\begin{thm}\label{proaddcla} Let $\beta$ be an interval-valued number and $X,Y\sqsubset U$, the the following statements hold.

(1) $\overline{SA}_{-1}^{\beta}(\emptyset)=\emptyset$, $\overline{SA}_{-1}^{\beta}(U)=U$.

(2) $\overline{SA}_{+1}^{\beta}(\emptyset)=\emptyset$, $\overline{SA}_{+1}^{\beta}(U)=U$.

(3) If $X\sqsubset Y$, then $\overline{SA}_{-1}^{\beta}(X)\sqsubset\overline{SA}_{-1}^{\beta}(Y)$ and $\overline{SA}_{+1}^{\beta}(X)\sqsubset\overline{SA}_{+1}^{\beta}(Y)$

(4) $\overline{SA}_{-1}^{\beta}(X)\sqcup\overline{SA}_{-1}^{\beta}(Y)\sqsubset\overline{SA}_{-1}^{\beta}(X\sqcup Y)$, $\overline{SA}_{+1}^{\beta}(X\sqcap Y)\sqsubset\overline{SA}_{+1}^{\beta}(X)\sqcap\overline{SA}_{+1}^{\beta}(Y)$.

(5) $\overline{SA}_{-1}^{\beta}(X)\sqcap\overline{SA}_{-1}^{\beta}(Y)=\overline{SA}_{-1}^{\beta}(X\sqcap Y)$, $\overline{SA}_{+1}^{\beta}(X\sqcup Y)=\overline{SA}_{+1}^{\beta}(X)\sqcup\overline{SA}_{+1}^{\beta}(Y)$.

(6) $\overline{SA}_{-1}^{\beta}(X^c)=(\overline{SA}_{+1}^{\beta}(X))^c$, $\overline{SA}_{+1}^{\beta}(X^c)=(\overline{SA}_{-1}^{\beta}(X))^c$.

(7) $\overline{SA}_{-1}^{\beta}(X)\sqsubset X\sqsubset\overline{SA}_{+1}^{\beta}(X)$.
\end{thm}

{\bf\slshape Proof} (1) and (2) are obvious.

(3) If $x\in\overline{SA}_{-1}^{\beta}(X)$, then $\overline{SN}^{\beta}_x\sqsubset X\sqsubset Y$, i.e., $x\in\overline{SA}_{-1}^{\beta}(Y)$. Then $\overline{SA}_{-1}^{\beta}(X)\sqsubset \overline{SA}_{-1}^{\beta}(Y)$.

If $x\in\overline{SA}_{+1}^{\beta}(X)$, then $\emptyset\neq\overline{SN}^{\beta}_x\sqcap X\sqsubset \overline{SN}^{\beta}_x\sqcap Y$, i.e., $x\in\overline{SA}_{+1}^{\beta}(Y)$. Then $\overline{SA}_{+1}^{\beta}(X)\sqsubset \overline{SA}_{+1}^{\beta}(Y)$.

(4) Since $X\sqsubset X\sqcup Y$ and $Y\sqsubset X\sqcup Y$, by the results of (3), $\overline{SA}_{-1}^{\beta}(X)\sqsubset\overline{SA}_{-1}^{\beta}(X\sqcup Y)$ and $\overline{SA}_{-1}^{\beta}(Y)\sqsubset\overline{SA}_{-1}^{\beta}(X\sqcup Y)$. Then  $\overline{SA}_{-1}^{\beta}(X)\sqcup\overline{SA}_{-1}^{\beta}(Y)\sqsubset\overline{SA}_{-1}^{\beta}(X\sqcup Y)$.

Since $X\sqcap Y\sqsubset X$ and $X\sqcap Y\sqsubset Y$, by the results of (3), $\overline{SA}_{+1}^{\beta}(X\sqcap Y)\sqsubset\overline{SA}_{+1}^{\beta}(X)$ and $\overline{SA}_{+1}^{\beta}(X\sqcap Y)\sqsubset\overline{SA}_{+1}^{\beta}(Y)$. Then $\overline{SA}_{+1}^{\beta}(X\sqcap Y)\sqsubset\overline{SA}_{+1}^{\beta}(X)\sqcap\overline{SA}_{+1}^{\beta}(Y)$.

(5) If $x\in\overline{SA}_{-1}^{\beta}(X)\sqcap\overline{SA}_{-1}^{\beta}(Y)$, i.e., $x\in\overline{SA}_{-1}^{\beta}(X)$ and $x\in\overline{SA}_{-1}^{\beta}(Y)$, then $\overline{SN}^{\beta}_x\sqsubset X$ and $\overline{SN}^{\beta}_x\sqsubset Y$. It means that $\overline{SN}^{\beta}_x\sqsubset X\sqcap Y$, i.e., $x\in\overline{SA}_{-1}^{\beta}(X\sqcap Y)$. Then $\overline{SA}_{-1}^{\beta}(X)\sqcap\overline{SA}_{-1}^{\beta}(Y)\sqsubset\overline{SA}_{-1}^{\beta}(X\sqcap Y)$.

On the other hand, $X\sqcap Y\sqsubset X$ and $X\sqcap Y\sqsubset Y$, by the results of (3), $\overline{SA}_{-1}^{\beta}(X\sqcap Y)\sqsubset\overline{SA}_{-1}^{\beta}(X)\sqcap\overline{SA}_{-1}^{\beta}(Y)$.

Hence, $\overline{SA}_{-1}^{\beta}(X\sqcap Y)=\overline{SA}_{-1}^{\beta}(X)\sqcap\overline{SA}_{-1}^{\beta}(Y)$.

If $x\in\overline{SA}_{+1}^{\beta}(X\sqcup Y)$, i.e., $\overline{SN}^{\beta}_x\sqcap(X\sqcup Y)\neq\emptyset$, then there is $y\in U$ such that $y\in\overline{SN}^{\beta}_x$ and $y\in X\sqcup Y$, it means that at least one of $y\in X$ and $y\in Y$ holds, then we can obtain that at least one of $\overline{SN}^{\beta}_x\sqcap X\neq\emptyset$ and $\overline{SN}^{\beta}_x\sqcap Y\neq\emptyset$ holds, i.e., $x\in\overline{SA}_{+1}^{\beta}(X)\sqcup\overline{SA}_{+1}^{\beta}(Y)$. Then $\overline{SA}_{+1}^{\beta}(X\sqcup Y)\sqsubset\overline{SA}_{+1}^{\beta}(X)\sqcup\overline{SA}_{+1}^{\beta}(Y)$.

On the other hand, $X\sqsubset X\sqcup Y$ and $Y\sqsubset X\sqcup Y$, by the results of (3), $\overline{SA}_{+1}^{\beta}(X)\sqcup\overline{SA}_{+1}^{\beta}(Y)\sqsubset\overline{SA}_{+1}^{\beta}(X\sqcup Y)$.

Hence, $\overline{SA}_{+1}^{\beta}(X)\sqcup\overline{SA}_{+1}^{\beta}(Y)=\overline{SA}_{+1}^{\beta}(X\sqcup Y)$.

(6) If $x\in\overline{SA}_{-1}^{\beta}(X^c)$, i.e., $\overline{SN}^{\beta}_x\sqsubset X^c$, then $\overline{SN}^{\beta}_x\sqcap X=\emptyset$, i.e., $x\not\in\overline{SA}_{+1}^{\beta}(X)$, thus $x\in(\overline{SA}_{+1}^{\beta}(X))^c$. Then  $\overline{SA}_{-1}^{\beta}(X^c)\sqsubset (\overline{SA}_{+1}^{\beta}(X))^c$.

On the other hand, $x\in(\overline{SA}_{+1}^{\beta}(X))^c$, i.e., $x\not\in\overline{SA}_{+1}^{\beta}(X)$, thus $\overline{SN}^{\beta}_x\sqcap X=\emptyset$. Then $\overline{SN}^{\beta}_x\sqsubset X^c$, i.e., $x\in\overline{SA}_{-1}^{\beta}(X^c)$. Then  $(\overline{SA}_{+1}^{\beta}(X))^c\sqsubset\overline{SA}_{-1}^{\beta}(X^c)$.

Hence, $\overline{SA}_{-1}^{\beta}(X^c)=(\overline{SA}_{+1}^{\beta}(X))^c$.

If  $x\in\overline{SA}_{+1}^{\beta}(X^c)$, i.e., $\overline{SN}^{\beta}_x\sqcap X^c\neq\emptyset$, then $\overline{SN}^{\beta}_x\not\sqsubset X$, i.e., $x\not\in\overline{SA}_{-1}^{\beta}(X)$ and $x\in(\overline{SA}_{-1}^{\beta}(X))^c$. Then $\overline{SA}_{+1}^{\beta}(X^c)\sqsubset(\overline{SA}_{-1}^{\beta}(X))^c$.

On the other hand, if $x\in(\overline{SA}_{-1}^{\beta}(X))^c$, i.e., $x\not\in\overline{SA}_{-1}^{\beta}(X)$ and $\overline{SN}^{\beta}_x\not\sqsubset X$. Thus  $\overline{SN}^{\beta}_x\sqcap X^c\neq\emptyset$, i.e., $x\in\overline{SA}_{+1}^{\beta}(X^c)$. Then $(\overline{SA}_{-1}^{\beta}(X))^c\sqsubset\overline{SA}_{+1}^{\beta}(X^c)$.

Hence, $\overline{SA}_{+1}^{\beta}(X^c)=(\overline{SA}_{-1}^{\beta}(X))^c$.

(7) By Theorem \ref{pro320}, $x\in\overline{SN}^{\beta}_x$ for all $x\in U$.

For a random $y\not\in X$ ($y\in X^c$), since $y\in\overline{SN}^{\beta}_y$, then $\overline{SN}^{\beta}_y\not\sqsubset X$, i.e., $y\not\in\overline{SA}_{-1}^{\beta}(X)$ ($y\in(\overline{SA}_{-1}^{\beta}(X))^c$), then $X^c\sqsubset(\overline{SA}_{-1}^{\beta}(X))^c$, i.e.,  $\overline{SA}_{-1}^{\beta}(X)\sqsubset X$.

For all $x\in X$, since $x\in\overline{SN}^{\beta}_x$, then $\emptyset\neq\{x\}\sqsubset\overline{SN}^{\beta}_x\sqcap X$, i.e., $x\in\overline{SA}_{+1}^{\beta}(X)$. Then $X\sqsubset\overline{SA}_{+1}^{\beta}(X)$.

To sum up, $\overline{SA}_{-1}^{\beta}(X)\sqsubset X\sqsubset\overline{SA}_{+1}^{\beta}(X)$.$\blacksquare$

\begin{thm} Let $\beta$ be an interval-valued number, $X\sqsubset U$, $A\sqsubset E$, $B\sqsubset E$.
In $T_{1}=(U,F,A)_{\beta}$ and $T_{2}=(U,F,B)_{\beta}$, if $\overline{SN}^{1\beta}_x=\overline{SN}^{2\beta}_x$ for all $x\in U$, then $\overline{SA}_{-1}^{1\beta}(X)=\overline{SA}_{-1}^{2\beta}(X)$ and $\overline{SA}_{+1}^{1\beta}(X)=\overline{SA}_{+1}^{2\beta}(X)$.

{\bf\slshape Proof}  (1) Since $\overline{SN}^{1\beta}_x=\overline{SN}^{2\beta}_x$ for all $x\in U$, then
$\overline{SN}^{1\beta}_x\sqsubset X$ if and only if $\overline{SN}^{2\beta}_x\sqsubset X$, then $\overline{SA}_{-1}^{1\beta}(X)=\overline{SA}_{-1}^{2\beta}(X)$.

$\overline{SN}^{1\beta}_x\sqcap X\neq\emptyset$ if and only if $\overline{SN}^{2\beta}_x\sqcap X\neq\emptyset$, then $\overline{SA}_{+1}^{1\beta}(X)=\overline{SA}_{+1}^{2\beta}(X)$.$\blacksquare$
\end{thm}

\subsection{The 2nd kind of interval-valued fuzzy soft $\beta$-coverings based fuzzy rough sets}

\begin{defn}\label{defsm} Let $\beta$ be an interval-valued number, and $U$ and $E$ be the sets of objects and parameters, respectively.
$\widetilde{SM}^{\beta}_x$ is called the interval-valued fuzzy soft complementary  $\beta$-neighborhood of $x$ in $(U,F,A)_{\beta}$, where $\widetilde{SM}^{\beta}_x(y)=\widetilde{SN}^{\beta}_y(x)$.

\end{defn}

\begin{defn}  Let $\beta$ be an interval-valued number, $X\in IF(U)$, $x,y\in U$. $\widetilde{SA}_{-2}^{\beta}(X)$ and $\widetilde{SA}_{+2}^{\beta}(X)$ are obtained by
\begin{equation*}
\left\{
\begin{array}{l}
\widetilde{SA}_{-2}^{\beta}(X)(x)=\bigwedge\limits_{y\in U}\{(\widetilde{SM}_x^{\beta})^c(y)\vee X(y)\}, \\
\widetilde{SA}_{+2}^{\beta}(X)(x)=\bigvee\limits_{y\in U}\{\widetilde{SM}_x^{\beta}(y)\wedge X(y)\}.
\end{array}
\right.
\end{equation*}

If $\widetilde{SA}_{-2}^{\beta}(X)\neq\widetilde{SA}_{+2}^{\beta}(X)$, then $X$ is called the 2nd type of interval-valued fuzzy
soft $\beta$-covering based interval-valued fuzzy rough sets; otherwise $X$ is called
interval-valued fuzzy definable.

\end{defn}

\begin{thm}\label{} Let $\beta$ be an interval-valued number, $X,Y\in IF(U)$, $x,y\in U$,  the following statements hold.

(1)  $\widetilde{SA}_{-2}^{\beta}(I^U)=I^U$, $\widetilde{SA}_{+2}^{\beta}(I^{\emptyset})= I^{\emptyset}$.

(2) $\widetilde{SA}_{-2}^{\beta}(X^c)=(\widetilde{SA}_{+2}^{\beta}(X))^c$, $\widetilde{SA}_{+2}^{\beta}(X^c)=(\widetilde{SA}_{-2}^{\beta}(X))^c$.

(3)  $\widetilde{SA}_{-2}^{\beta}(X\cap Y)=\widetilde{SA}_{-2}^{\beta}(X)\cap \widetilde{SA}_{-2}^{\beta}(Y)$, $\widetilde{SA}_{+2}^{\beta}(X\cup Y)=\widetilde{SA}_{+2}^{\beta}(X)\cup \widetilde{SA}_{+2}^{\beta}(Y)$.

(4) If $X\subset Y$, then $\widetilde{SA}_{-2}^{\beta}(X)\subset\widetilde{SA}_{-2}^{\beta}(Y)$, $\widetilde{SA}_{+2}^{\beta}(X)\subset\widetilde{SA}_{+2}^{\beta}(Y)$.

(5) $\widetilde{SA}_{-2}^{\beta}(X)\cup\widetilde{SA}_{-2}^{\beta}(Y)\subset\widetilde{SA}_{-2}^{\beta}(X\cup Y)$, $\widetilde{SA}_{+2}^{\beta}(X\cap Y)\subset\widetilde{SA}_{+2}^{\beta}(X)\cap\widetilde{SA}_{+2}^{\beta}(Y)$.

(6) If $(\widetilde{SM}_x^{\beta})^c(x)\leqslant X(x)\leqslant\widetilde{SM}_x^{\beta}(x)$ for all $x\in U$, then $\widetilde{SA}_{-2}^{\beta}(X)\subset X\subset\widetilde{SA}_{+2}^{\beta}(X)$.

(7) If $(\widetilde{SM}_x^{\beta})^c(x)\leqslant X(x)\leqslant\widetilde{SM}_x^{\beta}(x)$ for all $x\in U$, then $\widetilde{SA}_{-2}^{\beta}(\widetilde{SA}_{-2}^{\beta}(X))\subset\widetilde{SA}_{-2}^{\beta}(X)\subset X\subset\widetilde{SA}_{+2}^{\beta}(X)\subset\widetilde{SA}_{+2}^{\beta}(\widetilde{SA}_{+2}^{\beta}(X))$.

(8) If $X\subset Y$, $\widetilde{SA}_{-2}^{\beta}(X)\cup\widetilde{SA}_{-2}^{\beta}(Y)=\widetilde{SA}_{-2}^{\beta}(X\cup Y)$, $\widetilde{SA}_{+2}^{\beta}(X\cap Y)=\widetilde{SA}_{+2}^{\beta}(X)\cap\widetilde{SA}_{+2}^{\beta}(Y)$.

{\bf\slshape Proof} It can be proved similarly to Theorem \ref{proada} by using  $\widetilde{SM}^{\beta}_x$ and  $(\widetilde{SM}_x^{\beta})^c$ to instead of $\widetilde{SN}^{\beta}_x$ and $(\widetilde{SN}^{\beta}_x)^c$ in  Theorem \ref{proada}, respectively.$\blacksquare$
\end{thm}

\begin{defn}\label{defcla2}  Let $\beta$ be an interval-valued number and $X\sqsubset U$ be an object set, $x,y\in U$.
$\overline{SA}_{-2}^{\beta}(X)=\{x\in U:\overline{SM}^{\beta}_x\sqsubset X\}$ and $\overline{SA}_{+2}^{\beta}(X)=\{x\in U:\overline{SM}^{\beta}_x\sqcap X\neq\emptyset\}$ are the lower approximation and upper approximation of the object set $X$ in $(U,F,A)_{\beta}$, respectively.

If $\overline{SA}_{-2}^{\beta}(X)\neq \overline{SA}_{+2}^{\beta}(X)$, then $X$ is called the 2nd type of soft
$\beta$-covering based rough sets. Otherwise, $X$ is called definable.
\end{defn}

\subsection{The 3rd kind of interval-valued fuzzy soft $\beta$-coverings based fuzzy rough sets}

\begin{defn}  Let $\beta$ be an interval-valued number, $X\in IF(U)$, $x,y\in U$. $\widetilde{SA}_{-3}^{\beta}(X)$ and $\widetilde{SA}_{+3}^{\beta}(X)$ are obtained by
\begin{equation*}
\left\{
\begin{array}{l}
\widetilde{SA}_{-3}^{\beta}(X)(x)=\bigwedge\limits_{y\in U}\{(\widetilde{SN}_x^{\beta})^c(y)\vee(\widetilde{SM}_x^{\beta})^c(y)\vee X(y)\}, \\
\widetilde{SA}_{+3}^{\beta}(X)(x)=\bigvee\limits_{y\in U}\{\widetilde{SN}_x^{\beta}(y)\wedge\widetilde{SM}_x^{\beta}(y)\wedge X(y)\}.
\end{array}
\right.
\end{equation*}

\end{defn}

\begin{thm}\label{} Let $\beta$ be an interval-valued number, $X,Y\in IF(U)$, $x,y\in U$,  the following statements hold.

(1)  $\widetilde{SA}_{-3}^{\beta}(I^U)=I^U$, $\widetilde{SA}_{+3}^{\beta}(I^{\emptyset})= I^{\emptyset}$.

(2) $\widetilde{SA}_{-3}^{\beta}(X^c)=(\widetilde{SA}_{+3}^{\beta}(X))^c$, $\widetilde{SA}_{+3}^{\beta}(X^c)=(\widetilde{SA}_{-3}^{\beta}(X))^c$.

(3)  $\widetilde{SA}_{-3}^{\beta}(X\cap Y)=\widetilde{SA}_{-3}^{\beta}(X)\cap \widetilde{SA}_{-3}^{\beta}(Y)$, $\widetilde{SA}_{+3}^{\beta}(X\cup Y)=\widetilde{SA}_{+3}^{\beta}(X)\cup \widetilde{SA}_{+3}^{\beta}(Y)$.

(4) If $X\subset Y$, then $\widetilde{SA}_{-3}^{\beta}(X)\subset\widetilde{SA}_{-3}^{\beta}(Y)$, $\widetilde{SA}_{+3}^{\beta}(X)\subset\widetilde{SA}_{+3}^{\beta}(Y)$.

(5) $\widetilde{SA}_{-3}^{\beta}(X)\cup\widetilde{SA}_{-3}^{\beta}(Y)\subset\widetilde{SA}_{-3}^{\beta}(X\cup Y)$, $\widetilde{SA}_{+3}^{\beta}(X\cap Y)\subset\widetilde{SA}_{+3}^{\beta}(X)\cap\widetilde{SA}_{+3}^{\beta}(Y)$.

(6) If $(\widetilde{SN}_x^{\beta}\cap\widetilde{SM}_x^{\beta})^c(x)\leqslant X(x)\leqslant(\widetilde{SN}_x^{\beta}\cap\widetilde{SM}_x^{\beta})(x)$ for all $x\in U$, then $\widetilde{SA}_{-3}^{\beta}(X)\subset X\subset\widetilde{SA}_{+3}^{\beta}(X)$.

(7) If $(\widetilde{SN}_x^{\beta}\cap\widetilde{SM}_x^{\beta})^c(x)\leqslant X(x)\leqslant(\widetilde{SN}_x^{\beta}\cap\widetilde{SM}_x^{\beta})(x)$ for all $x\in U$, then $\widetilde{SA}_{-3}^{\beta}(\widetilde{SA}_{-3}^{\beta}(X))\subset\widetilde{SA}_{-3}^{\beta}(X)\subset X\subset\widetilde{SA}_{+3}^{\beta}(X)\subset\widetilde{SA}_{+3}^{\beta}(\widetilde{SA}_{+3}^{\beta}(X))$.

(8) If $X\subset Y$, $\widetilde{SA}_{-3}^{\beta}(X)\cup\widetilde{SA}_{-3}^{\beta}(Y)=\widetilde{SA}_{-3}^{\beta}(X\cup Y)$, $\widetilde{SA}_{+3}^{\beta}(X\cap Y)=\widetilde{SA}_{+3}^{\beta}(X)\cap\widetilde{SA}_{+3}^{\beta}(Y)$.

{\bf\slshape Proof}  It can be proved similarly to Theorem \ref{proada} by using  $\widetilde{SN}_x^{\beta}\cap\widetilde{SM}_x^{\beta}$ and  $(\widetilde{SN}_x^{\beta}\cap\widetilde{SM}_x^{\beta})^c$ to instead of $\widetilde{SN}^{\beta}_x$ and $(\widetilde{SN}^{\beta}_x)^c$ in  Theorem \ref{proada}, respectively.$\blacksquare$
\end{thm}

\begin{defn}\label{defcla3}  Let $\beta$ be an interval-valued number and $X\sqsubset U$ be an object set, $x,y\in U$.
$\overline{SA}_{-3}^{\beta}(X)=\{x\in U:\overline{SN}^{\beta}_x\sqsubset X\ or\ \overline{SM}^{\beta}_x\sqsubset X\}$ and $\overline{SA}_{+3}^{\beta}(X)=\{x\in U:\overline{SN}^{\beta}_x\sqcap X\neq\emptyset\ and\ \overline{SM}^{\beta}_x\sqcap X\neq\emptyset\}$ are the lower approximation and upper approximation of the object set $X$ in $(U,F,A)_{\beta}$, respectively.

If $\overline{SA}_{-3}^{\beta}(X)\neq \overline{SA}_{+3}^{\beta}(X)$, then $X$ is called the 3rd type of soft
$\beta$-covering based rough sets. Otherwise, $X$ is called definable.
\end{defn}

\subsection{The 4th kind of interval-valued fuzzy soft $\beta$-coverings based fuzzy rough sets}

\begin{defn}  Let $\beta$ be an interval-valued number, $X\in IF(U)$, $x,y\in U$. $\widetilde{SA}_{-4}^{\beta}(X)$ and $\widetilde{SA}_{+4}^{\beta}(X)$ are obtained by
\begin{equation*}
\left\{
\begin{array}{l}
\widetilde{SA}_{-4}^{\beta}(X)(x)=\bigwedge\limits_{y\in U}\{(\widetilde{SN}_x^{\beta})^c(y)\wedge(\widetilde{SM}_x^{\beta})^c(y)\vee X(y)\}, \\
\widetilde{SA}_{+4}^{\beta}(X)(x)=\bigvee\limits_{y\in U}\{\widetilde{SN}_x^{\beta}(y)\vee\widetilde{SM}_x^{\beta}(y)\wedge X(y)\}.
\end{array}
\right.
\end{equation*}

\end{defn}

\begin{thm}\label{} Let $\beta$ be an interval-valued number, $X,Y\in IF(U)$, $x,y\in U$,  the following statements hold.

(1)  $\widetilde{SA}_{-4}^{\beta}(I^U)=I^U$, $\widetilde{SA}_{+4}^{\beta}(I^{\emptyset})= I^{\emptyset}$.

(2) $\widetilde{SA}_{-4}^{\beta}(X^c)=(\widetilde{SA}_{+4}^{\beta}(X))^c$, $\widetilde{SA}_{+4}^{\beta}(X^c)=(\widetilde{SA}_{-4}^{\beta}(X))^c$.

(3)  $\widetilde{SA}_{-4}^{\beta}(X\cap Y)=\widetilde{SA}_{-4}^{\beta}(X)\cap \widetilde{SA}_{-4}^{\beta}(Y)$, $\widetilde{SA}_{+4}^{\beta}(X\cup Y)=\widetilde{SA}_{+4}^{\beta}(X)\cup \widetilde{SA}_{+4}^{\beta}(Y)$.

(4) If $X\subset Y$, then $\widetilde{SA}_{-4}^{\beta}(X)\subset\widetilde{SA}_{-4}^{\beta}(Y)$, $\widetilde{SA}_{+4}^{\beta}(X)\subset\widetilde{SA}_{+4}^{\beta}(Y)$.

(5) $\widetilde{SA}_{-4}^{\beta}(X)\cup\widetilde{SA}_{-4}^{\beta}(Y)\subset\widetilde{SA}_{-4}^{\beta}(X\cup Y)$, $\widetilde{SA}_{+4}^{\beta}(X\cap Y)\subset\widetilde{SA}_{+4}^{\beta}(X)\cap\widetilde{SA}_{+4}^{\beta}(Y)$.

(6) If $(\widetilde{SN}_x^{\beta}\cup\widetilde{SM}_x^{\beta})^c(x)\leqslant X(x)\leqslant(\widetilde{SN}_x^{\beta}\cup\widetilde{SM}_x^{\beta})(x)$ for all $x\in U$, then $\widetilde{SA}_{-4}^{\beta}(X)\subset X\subset\widetilde{SA}_{+4}^{\beta}(X)$.

(7) If $(\widetilde{SN}_x^{\beta}\cup\widetilde{SM}_x^{\beta})^c(x)\leqslant X(x)\leqslant(\widetilde{SN}_x^{\beta}\cup\widetilde{SM}_x^{\beta})(x)$ for all $x\in U$, then $\widetilde{SA}_{-4}^{\beta}(\widetilde{SA}_{-4}^{\beta}(X))\subset\widetilde{SA}_{-4}^{\beta}(X)\subset X\subset\widetilde{SA}_{+4}^{\beta}(X)\subset\widetilde{SA}_{+4}^{\beta}(\widetilde{SA}_{+4}^{\beta}(X))$.

(8) If $X\subset Y$, $\widetilde{SA}_{-4}^{\beta}(X)\cup\widetilde{SA}_{-4}^{\beta}(Y)=\widetilde{SA}_{-4}^{\beta}(X\cup Y)$, $\widetilde{SA}_{+4}^{\beta}(X\cap Y)=\widetilde{SA}_{+4}^{\beta}(X)\cap\widetilde{SA}_{+4}^{\beta}(Y)$.

{\bf\slshape Proof}  It can be proved similarly to Theorem \ref{proada} by using  $\widetilde{SN}_x^{\beta}\cup\widetilde{SM}_x^{\beta}$ and  $(\widetilde{SN}_x^{\beta}\cup\widetilde{SM}_x^{\beta})^c$ to instead of $\widetilde{SN}^{\beta}_x$ and $(\widetilde{SN}^{\beta}_x)^c$ in  Theorem \ref{proada}, respectively.$\blacksquare$
\end{thm}

\begin{defn}\label{defcla4}  Let $\beta$ be an interval-valued number and $X\sqsubset U$ be an object set, $x,y\in U$.
$\overline{SA}_{-4}^{\beta}(X)=\{x\in U:\overline{SN}^{\beta}_x\sqsubset X\ and\ \overline{SM}^{\beta}_x\sqsubset X\}$ and $\overline{SA}_{+4}^{\beta}(X)=\{x\in U:\overline{SN}^{\beta}_x\sqcap X\neq\emptyset\ or\ \overline{SM}^{\beta}_x\sqcap X\neq\emptyset\}$ are the lower approximation and upper approximation of the object set $X$ in $(U,F,A)_{\beta}$, respectively.

If $\overline{SA}_{-4}^{\beta}(X)\neq \overline{SA}_{+4}^{\beta}(X)$, then $X$ is called the 4th type of soft
$\beta$-covering based rough sets. Otherwise, $X$ is called definable.
\end{defn}

\subsection{The relationships of four kinds of interval-valued fuzzy soft $\beta$-coverings based fuzzy rough sets}

\begin{thm}\label{prorelfu} Let $\beta$ be an interval-valued number, $X\in IF(U)$,   the following statements hold in $(U,F,A)_{\beta}$.

(1) $\widetilde{SA}_{-3}^{\beta}(X)=\widetilde{SA}_{-1}^{\beta}(X)\cup \widetilde{SA}_{-2}^{\beta}(X)$.

(2) $\widetilde{SA}_{+3}^{\beta}(X)=\widetilde{SA}_{+1}^{\beta}(X)\cap \widetilde{SA}_{+2}^{\beta}(X)$.

(3) $\widetilde{SA}_{-4}^{\beta}(X)=\widetilde{SA}_{-1}^{\beta}(X)\cap \widetilde{SA}_{-2}^{\beta}(X)$.

(4) $\widetilde{SA}_{+4}^{\beta}(X)=\widetilde{SA}_{+1}^{\beta}(X)\cup \widetilde{SA}_{+2}^{\beta}(X)$.

(5) $\widetilde{SA}_{-4}^{\beta}(X)\subset\widetilde{SA}_{-1}^{\beta}(X)\subset \widetilde{SA}_{-3}^{\beta}(X)$.

(6) $\widetilde{SA}_{-4}^{\beta}(X)\subset\widetilde{SA}_{-2}^{\beta}(X)\subset \widetilde{SA}_{-3}^{\beta}(X)$.

(7) $\widetilde{SA}_{+3}^{\beta}(X)\subset\widetilde{SA}_{+1}^{\beta}(X)\subset \widetilde{SA}_{+4}^{\beta}(X)$.

(8) $\widetilde{SA}_{+3}^{\beta}(X)\subset\widetilde{SA}_{+2}^{\beta}(X)\subset \widetilde{SA}_{+4}^{\beta}(X)$.

{\bf\slshape Proof} (1) For all $x,y\in U$, $\widetilde{SA}_{-1}^{\beta}(X)(x)\vee \widetilde{SA}_{-2}^{\beta}(X)(x)=(\bigwedge\limits_{y\in U}\{(\widetilde{SN}_x^{\beta})^c(y)\vee X(y)\})\vee(\bigwedge\limits_{y\in U}\{(\widetilde{SM}_x^{\beta})^c(y)\vee X(y)\})=\bigwedge\limits_{y\in U}\{(\widetilde{SN}_x^{\beta})^c(y)\vee X(y)\vee(\widetilde{SM}_x^{\beta})^c(y)\vee X(y)\}=\bigwedge\limits_{y\in U}\{(\widetilde{SN}_x^{\beta})^c(y)\vee(\widetilde{SM}_x^{\beta})^c(y)\vee X(y)\}=\widetilde{SA}_{-3}^{\beta}(X)(x)$.

(2) For all $x,y\in U$, $\widetilde{SA}_{+1}^{\beta}(X)(x)\wedge \widetilde{SA}_{+2}^{\beta}(X)(x)=(\bigvee\limits_{y\in U}\{\widetilde{SN}_x^{\beta}(y)\wedge X(y)\})\wedge(\bigvee\limits_{y\in U}\{\widetilde{SM}_x^{\beta}(y)\wedge X(y)\})=\bigvee\limits_{y\in U}\{\widetilde{SN}_x^{\beta}(y)\wedge\widetilde{SM}_x^{\beta}(y)\wedge X(y)\wedge X(y)\}=\widetilde{SA}_{+3}^{\beta}(X)(x)$.

(3) For all $x,y\in U$, $\widetilde{SA}_{-1}^{\beta}(X)(x)\wedge\widetilde{SA}_{-2}^{\beta}(X)(x)=(\bigwedge\limits_{y\in U}\{(\widetilde{SN}_x^{\beta})^c(y)\vee X(y)\})\wedge(\bigwedge\limits_{y\in U}\{(\widetilde{SM}_x^{\beta})^c(y)\vee X(y)\})=\bigwedge\limits_{y\in U}\{[(\widetilde{SM}_x^{\beta})^c(y)\vee X(y)]\wedge[(\widetilde{SN}_x^{\beta})^c(y)\vee X(y)]\}=\bigwedge\limits_{y\in U}\{[(\widetilde{SM}_x^{\beta})^c(y)\wedge(\widetilde{SN}_x^{\beta})^c(y)]\vee X(y)\}=\widetilde{SA}_{-4}^{\beta}(X)(x)$.

(4) For all $x,y\in U$, $\widetilde{SA}_{+1}^{\beta}(X)(x)\vee\widetilde{SA}_{+2}^{\beta}(X)(x)=(\bigvee\limits_{y\in U}\{\widetilde{SN}_x^{\beta}(y)\wedge X(y)\})\vee(\bigvee\limits_{y\in U}\{\widetilde{SM}_x^{\beta}(y)\wedge X(y)\})=\bigvee\limits_{y\in U}\{[\widetilde{SN}_x^{\beta}(y)\vee\widetilde{SM}_x^{\beta}(y)]\wedge X(y)\}=\widetilde{SA}_{+4}^{\beta}(X)(x)$.

(5) For all $y\in U$, $(\widetilde{SN}_x^{\beta})^c(y)\wedge(\widetilde{SM}_x^{\beta})^c(y)\vee X(y)\leqslant(\widetilde{SN}_x^{\beta})^c(y)\vee X(y)\leqslant(\widetilde{SN}_x^{\beta})^c(y)\vee(\widetilde{SM}_x^{\beta})^c(y)\vee X(y)$.

(6) For all $y\in U$, $(\widetilde{SN}_x^{\beta})^c(y)\wedge(\widetilde{SM}_x^{\beta})^c(y)\vee X(y)\leqslant(\widetilde{SM}_x^{\beta})^c(y)\vee X(y)\leqslant(\widetilde{SN}_x^{\beta})^c(y)\vee(\widetilde{SM}_x^{\beta})^c(y)\vee X(y)$.

(7) For all $y\in U$, $\widetilde{SN}_x^{\beta}(y)\wedge\widetilde{SM}_x^{\beta}(y)\wedge X(y)\leqslant\widetilde{SN}_x^{\beta}(y)\wedge X(y)\leqslant\widetilde{SN}_x^{\beta}(y)\vee\widetilde{SM}_x^{\beta}(y)\wedge X(y)$.

(8)  For all $y\in U$, $\widetilde{SN}_x^{\beta}(y)\wedge\widetilde{SM}_x^{\beta}(y)\wedge X(y)\leqslant\widetilde{SM}_x^{\beta}(y)\wedge X(y)\leqslant\widetilde{SN}_x^{\beta}(y)\vee\widetilde{SM}_x^{\beta}(y)\wedge X(y)$.$\blacksquare$

\end{thm}

\begin{thm}\label{prorelcla} Let $\beta$ be an interval-valued number, $X\in IF(U)$,   the following statements hold in $(U,F,A)_{\beta}$.

(1) $\overline{SA}_{-3}^{\beta}(X)=\overline{SA}_{-1}^{\beta}(X)\sqcup \overline{SA}_{-2}^{\beta}(X)$.

(2) $\overline{SA}_{+3}^{\beta}(X)=\overline{SA}_{+1}^{\beta}(X)\sqcap \overline{SA}_{+2}^{\beta}(X)$.

(3) $\overline{SA}_{-4}^{\beta}(X)=\overline{SA}_{-1}^{\beta}(X)\sqcap \overline{SA}_{-2}^{\beta}(X)$.

(4) $\overline{SA}_{+4}^{\beta}(X)=\overline{SA}_{+1}^{\beta}(X)\sqcup \overline{SA}_{+2}^{\beta}(X)$.

(5) $\overline{SA}_{-4}^{\beta}(X)\sqsubset\overline{SA}_{-1}^{\beta}(X)\sqsubset \overline{SA}_{-3}^{\beta}(X)$.

(6) $\overline{SA}_{-4}^{\beta}(X)\sqsubset\overline{SA}_{-2}^{\beta}(X)\sqsubset \overline{SA}_{-3}^{\beta}(X)$.

(7) $\overline{SA}_{+3}^{\beta}(X)\sqsubset\overline{SA}_{+1}^{\beta}(X)\sqsubset \overline{SA}_{+4}^{\beta}(X)$.

(8) $\overline{SA}_{+3}^{\beta}(X)\sqsubset\overline{SA}_{+2}^{\beta}(X)\sqsubset \overline{SA}_{+4}^{\beta}(X)$.

{\bf\slshape Proof}   It is clear from Definition \ref{defcla1}, \ref{defcla2}, \ref{defcla3} and \ref{defcla4}.$\blacksquare$

\end{thm}

\begin{thm}\label{} Let $\beta$ be an interval-valued number, $X\in IF(U)$, if $(\widetilde{SN}_x^{\beta})^c(x)\leqslant X(x)\leqslant \widetilde{SN}_x^{\beta}(x)$ for all $x\in U$, then
$\widetilde{SA}_{-4}^{\beta}(X)\subset\widetilde{SA}_{-2}^{\beta}(X)\subset\widetilde{SA}_{-3}^{\beta}(X)\subset X\subset\widetilde{SA}_{+3}^{\beta}(X)\subset\widetilde{SA}_{+1}^{\beta}(X)\subset\widetilde{SA}_{+4}^{\beta}(X)$.

{\bf\slshape Proof} For all $x,y\in U$, $\widetilde{SM}_x^{\beta}(y)=\widetilde{SN}_y^{\beta}(x)$, then $\widetilde{SM}_x^{\beta}(x)=\widetilde{SN}_x^{\beta}(x)$ and $(\widetilde{SM}_x^{\beta})^c(x)=(\widetilde{SN}_x^{\beta})^c(x)$.

For all $x\in U$, $(\widetilde{SN}_x^{\beta})^c(x)\leqslant X(x)$, then $(\widetilde{SN}_x^{\beta})^c(x)\vee(\widetilde{SM}_x^{\beta})^c(x)\vee X(x)=(\widetilde{SN}_x^{\beta})^c(x)\vee X(x)\leqslant X(x)$. Then, $\widetilde{SA}_{-3}^{\beta}(X)(x)=\bigwedge\limits_{y\in U}\{(\widetilde{SN}_x^{\beta})^c(y)\vee(\widetilde{SM}_x^{\beta})^c(y)\vee X(y)\}\leqslant (\widetilde{SN}_x^{\beta})^c(x)\vee(\widetilde{SM}_x^{\beta})^c(x)\vee X(x)\leqslant X(x)$ for all $x\in U$, i.e., $\widetilde{SA}_{-3}^{\beta}(X)\subset X$.
 
For all $x\in U$, $X(x)\leqslant \widetilde{SN}_x^{\beta}(x)$, then $X(x)\leqslant \widetilde{SN}_x^{\beta}(x)\wedge X(x)=\widetilde{SN}_x^{\beta}(x)\wedge \widetilde{SM}_x^{\beta}(x)\wedge X(x)$. Then, $\widetilde{SA}_{+3}^{\beta}(X)(x)=\bigvee\limits_{y\in U}\{\widetilde{SN}_x^{\beta}(y)\wedge\widetilde{SM}_x^{\beta}(y)\wedge X(y)\}\geqslant \widetilde{SN}_x^{\beta}(x)\wedge \widetilde{SM}_x^{\beta}(x)\wedge X(x)\geqslant X(x)$ for all $x\in U$, i.e., $\widetilde{SA}_{+3}^{\beta}(X)\supset X$. 

To combine the results of Theorem \ref{prorelfu} (6) and (7), $\widetilde{SA}_{-4}^{\beta}(X)\subset\widetilde{SA}_{-2}^{\beta}(X)\subset\widetilde{SA}_{-3}^{\beta}(X)\subset X\subset\widetilde{SA}_{+3}^{\beta}(X)\subset\widetilde{SA}_{+1}^{\beta}(X)\subset\widetilde{SA}_{+4}^{\beta}(X)$ holds.$\blacksquare$

\end{thm}

\section*{References}

\bibliographystyle{plain}
\bibliography{reference1}

\begin{thebibliography}{10}

\bibitem{albaity2023impact}
M.~Albaity, T.~Mahmood, and et~al.
\newblock Impact of machine learning and artificial intelligence in business
  based on intuitionistic fuzzy soft waspas method.
\newblock {\em Mathematics}, 11(6):1453, 2023.

\bibitem{alcantud2019n}
J.~C.~R. Alcantud, F.~Feng, and et~al.
\newblock An $n$-soft set approach to rough sets.
\newblock {\em IEEE Transactions on Fuzzy Systems}, 28(11):2996--3007, 2019.

\bibitem{atef2021three}
M.~Atef and S.~I. Nada.
\newblock On three types of soft fuzzy coverings based rough sets.
\newblock {\em Mathematics and Computers in Simulation}, 185:452--467, 2021.

\bibitem{zahedi}
Z.~K. Azadeh and K.~Adem.
\newblock Multi-attribute decision-making based on soft set theory: A
  systematic review.
\newblock {\em Soft Computing}, 23:6899--6920, 2019.

\bibitem{bayramov2023interval}
S.~Bayramov, {\c{C}}.~G. Aras, and et~al.
\newblock Interval-valued topology on soft sets.
\newblock {\em Axioms}, 12(7):692, 2023.

\bibitem{BJAT}
J.~B{\l}aszczy{\'n}ski, A.~T. de~Almeida~Filho, and et~al.
\newblock Auto loan fraud detection using dominance-based rough set approach
  versus machine learning methods.
\newblock {\em Expert Systems with Applications}, 163:113740, 2021.

\bibitem{bustince2013generation}
H.~Bustince, J.~Fernandez, and et~al.
\newblock Generation of linear orders for intervals by means of aggregation
  functions.
\newblock {\em Fuzzy Sets and Systems}, 220:69--77, 2013.

\bibitem{ccalik2021novel}
A.~{\c{C}}al{\i}k.
\newblock A novel pythagorean fuzzy ahp and fuzzy topsis methodology for green
  supplier selection in the industry 4.0 era.
\newblock {\em Soft Computing}, 25(3):2253--2265, 2021.

\bibitem{CHLT}
H.~Chen, T.~Li, and et~al.
\newblock A decision-theoretic rough set approach for dynamic data mining.
\newblock {\em IEEE Transactions on fuzzy Systems}, 23(6):1958--1970, 2015.

\bibitem{2023Intention}
K.~Chen, V.~Lou, and et~al.
\newblock Intention to use robotic exoskeletons by older people: A fuzzy-set
  qualitative comparative analysis approach.
\newblock {\em Computers in Human Behavior}, 141:107610, 2023.

\bibitem{cheng2023error}
H.~Cheng and L.~Gao.
\newblock Error exponent and strong converse for quantum soft covering.
\newblock {\em IEEE Transactions on Information Theory}, 2023.

\bibitem{dai2023interval}
J.~Dai, Z.~Wang, and et~al.
\newblock Interval-valued fuzzy discernibility pair approach for attribute
  reduction in incomplete interval-valued information systems.
\newblock {\em Information Sciences}, 642:119215, 2023.

\bibitem{fan2024overlap}
Y.~Fan, X.~Zhang, and et~al.
\newblock Overlap function-based fuzzy $\beta$-covering relations and fuzzy
  $\beta$-covering rough set models.
\newblock {\em International Journal of Approximate Reasoning}, 169:109164,
  2024.

\bibitem{feng2016soft}
F.~Feng, J.~Cho, and et~al.
\newblock Soft set based association rule mining.
\newblock {\em Knowledge-Based Systems}, 111:268--282, 2016.

\bibitem{feng2010application}
F.~Feng, Y.~Li, and et~al.
\newblock Application of level soft sets in decision making based on
  interval-valued fuzzy soft sets.
\newblock {\em Computers \& Mathematics with Applications}, 60(6):1756--1767,
  2010.

\bibitem{gorzalczany1987method}
M.~B. Gorza{\l}czany.
\newblock A method of inference in approximate reasoning based on
  interval-valued fuzzy sets.
\newblock {\em Fuzzy sets and systems}, 21(1):1--17, 1987.

\bibitem{2023Evaluation}
B.~Gueneri and M.~Deveci.
\newblock Evaluation of supplier selection in the defense industry using q-rung
  orthopair fuzzy set based edas approach.
\newblock {\em Expert Systems with Application}, 222:119846, 2023.

\bibitem{guo2022three}
D.~Guo, C.~Jiang, and et~al.
\newblock Three-way decision based on confidence level change in rough set.
\newblock {\em International Journal of Approximate Reasoning}, 143:57--77,
  2022.

\bibitem{jiang2013entropy}
Y.~Jiang, Y.~Tang, and et~al.
\newblock Entropy on intuitionistic fuzzy soft sets and on interval-valued
  fuzzy soft sets.
\newblock {\em Information Sciences}, 240:95--114, 2013.

\bibitem{KES}
E.~Kannan, S.~Ravikumar, and et~al.
\newblock Analyzing uncertainty in cardiotocogram data for the prediction of
  fetal risks based on machine learning techniques using rough set.
\newblock {\em Journal of Ambient Intelligence and Humanized Computing}, pages
  1--13, 2021.

\bibitem{lin2021risk}
S.~Lin, S.~Shen, and et~al.
\newblock Risk assessment and management of excavation system based on fuzzy
  set theory and machine learning methods.
\newblock {\em Automation in Construction}, 122:103490, 2021.

\bibitem{liu2021novel}
Y.~Liu, L.~Ma, and et~al.
\newblock A novel robust fuzzy mean-upm model for green closed-loop supply
  chain network design under distribution ambiguity.
\newblock {\em Applied Mathematical Modelling}, 92:99--135, 2021.

\bibitem{ma2016two}
L.~Ma.
\newblock Two fuzzy covering rough set models and their generalizations over
  fuzzy lattices.
\newblock {\em Fuzzy Sets and Systems}, 294:1--17, 2016.

\bibitem{movckovr2021approximations}
J.~Mo{\v{c}}ko{\v{r}} and P.~Hurt{\'\i}k.
\newblock Approximations of fuzzy soft sets by fuzzy soft relations with image
  processing application.
\newblock {\em Soft Computing}, 25(10):6915--6925, 2021.

\bibitem{MOLOD}
D.~Molodtsov.
\newblock Soft set theory-first results.
\newblock {\em Computers \& Mathematics With Applications}, 37:19--31, 1999.

\bibitem{mondal1999topology}
T.~K. Mondal and S.~K. Samanta.
\newblock Topology of interval-valued fuzzy sets.
\newblock {\em Indian Journal of Pure and Applied Mathematics}, 30:23--29,
  1999.

\bibitem{FF}
F.~Patras.
\newblock {\em The Essence of Numbers-Lecture Notes in Mathematics book
  series}.
\newblock Springer Press, 2020.

\bibitem{PALA82}
Z.~Pawlak.
\newblock Rough sets.
\newblock {\em International Journal of Applied Mathematics and Computer
  Science}, 11:341--356, 1982.

\bibitem{pkekala2021inclusion}
B.~P{\c{e}}kala, K.~Dyczkowski, and et~al.
\newblock Inclusion and similarity measures for interval-valued fuzzy sets
  based on aggregation and uncertainty assessment.
\newblock {\em Information Sciences}, 547:1182--1200, 2021.

\bibitem{peng2018algorithms}
X.~Peng and H.~Garg.
\newblock Algorithms for interval-valued fuzzy soft sets in emergency decision
  making based on wdba and codas with new information measure.
\newblock {\em Computers \& Industrial Engineering}, 119:439--452, 2018.

\bibitem{perez2015ordering}
R.~P{\'e}rez-Fern{\'a}ndez, P.~Alonso, and et~al.
\newblock Ordering finitely generated sets and finite interval-valued hesitant
  fuzzy sets.
\newblock {\em Information Sciences}, 325:375--392, 2015.

\bibitem{MTNA}
M.~Tishya and A.~Anitha.
\newblock Precipitation prediction by integrating rough set on fuzzy
  approximation space with deep learning techniques.
\newblock {\em Applied Soft Computing}, 139:110253, 2023.

\bibitem{WGLT}
G.~Wang, T.~Li, and et~al.
\newblock Double-local rough sets for efficient data mining.
\newblock {\em Information Sciences}, 571:475--498, 2021.

\bibitem{yang2016fuzzy}
B.~Yang and B.~Hu.
\newblock A fuzzy covering-based rough set model and its generalization over
  fuzzy lattice.
\newblock {\em Information Sciences}, 367:463--486, 2016.

\bibitem{yang2017some}
B.~Yang and B.~Hu.
\newblock On some types of fuzzy covering-based rough sets.
\newblock {\em Fuzzy sets and Systems}, 312:36--65, 2017.

\bibitem{yang2013}
Y.~Yang, X.~Tan, and et~al.
\newblock The multi-fuzzy soft set and its application in decision making.
\newblock {\em Applied Mathematical Modelling}, 37(7):4915--4923, 2013.

\bibitem{JYE}
J.~Ye, J.~Zhan, and Z.~Xu.
\newblock A novel multi-attribute decision-making method based on fuzzy rough
  sets.
\newblock {\em Computers \& Industrial Engineering}, 155:107136, 2021.

\bibitem{YUG}
G.~Yu.
\newblock Relationships between fuzzy approximation spaces and their
  uncertainty measures.
\newblock {\em Information Sciences}, 528:181--204, 2020.

\bibitem{zadeh1965fuzzy}
L.~A. Zadeh.
\newblock Fuzzy sets.
\newblock {\em Information and control}, 8(3):338--353, 1965.

\bibitem{zhang2009characterization}
H.~Zhang, W.~Zhang, and et~al.
\newblock On characterization of generalized interval-valued fuzzy rough sets
  on two universes of discourse.
\newblock {\em International Journal of Approximate Reasoning}, 51(1):56--70,
  2009.

\bibitem{ZKZJ}
K.~Zhang, J.~Zhan, and et~al.
\newblock On multicriteria decision-making method based on a fuzzy rough set
  model with fuzzy $\alpha$-neighborhoods.
\newblock {\em IEEE Transactions on Fuzzy Systems}, 29(9):2491--2505, 2020.

\bibitem{zhangzhan}
L.~Zhang and J.~Zhan.
\newblock Fuzzy soft $\beta$-covering based fuzzy rough sets and corresponding
  decision-making applications.
\newblock {\em International Journal of Machine Learning and Cybernetics},
  10(6):1487--1502, 2019.

\bibitem{zhang2020fuzzy}
X.~Zhang and J.~Wang.
\newblock Fuzzy $\beta$-covering approximation spaces.
\newblock {\em International Journal of Approximate Reasoning}, 126:27--47,
  2020.

\bibitem{zhang2024event}
Z.~Zhangu, X.~Yang, and et~al.
\newblock Event-triggered control for switched systems with sensor faults via
  adaptive fuzzy observer.
\newblock {\em Mathematics and Computers in Simulation}, 221:244--259, 2024.

\end{thebibliography}







\end{multicols}
\end{document}